\documentstyle[12pt,amscd,amssymb,epsf]{amsart}
	\newcommand{\zz}{{\Bbb Z}}
	\newcommand{\qq}{{\Bbb Q}}

\newtheorem{defn}{Definition}

\newtheorem{prop}[defn]{Proposition}

\begin{document}
\title{Integrality of the Averaged Jones Polynomial of Algebraically Split
Links}
\author{Hans U. Boden}
\address{Department of Mathematics and Statistics, McMaster University,
Hamilton, Ontario L8S 4K1 Canada}
\email{boden@@icarus.mcmaster.ca}
\maketitle

This note arose out of an attempt to prove a conjecture of Lin
and Wang concerning the the integrality of
the coefficients of the Taylor series expansion at $t=1$
of the averaged Jones polynomial
of algebraically split links.
This question comes up in the study of
Ohtsuki's invariants, $\lambda_n(M) \in \qq,$ defined for $M$ an
integral homology
3-sphere \cite{ohtsuki, lin-wang}.

\bigskip

Suppose that $L$ is an oriented link in $S^3$ with $\mu$ components.
Write the
averaged Jones polynomial of $L$ as a Taylor series at $t=1,$ i.e.
$$ \Phi(L;t) = \sum_{i=0}^\infty a_i (t-1)^i,$$
and set
$\phi_n(L)= (-2)^{\mu} a_{n+\mu}(L).$

Conjecture 4.1 of \cite{lin-wang} states that for $L$ an algebraically
split link (ASL),
$$n! \, \phi_n(L) \in 6 \zz.$$
This conjecture is verified for $n=1,2$ in \cite{lin-wang},
and we consider the case $n \ge 3$ here. We first establish that
$a_n(L) \in \zz$
whenever $L$ is a geometrically split link (GSL),
implying that $\phi_n(L) \in 2^{\mu} \zz,$ which is a priori stronger
than the conjecture in this case.
Nevertheless, Conjecture 4.1
is not true for ASLs.
The problem is the presence of additional factors of 2 in the denominator of
$\phi_n(L).$

The goal of this paper is to present two results, Proposition 1 for GSLs and
Proposition 2 for ASL, both giving integrality results for the coefficients
$a_n(L)$
which are sharp, as shown by examples (I) and (II).

\bigskip

For the definition of $\Phi(L;t)$, see pp.10--11 of \cite{lin-wang}.
It is roughly a sum of (normalized) Jones polynomials, summed over
{\it all} sublinks of $L,$ and it satisfies:

\medskip

\begin{enumerate}
\item[(i)] If $L$ is trivial or empty, then $\Phi(L;t)=1.$
\item[(ii)] If $L=L_1 \cup L_2$ is a geometric splitting,
then $\Phi(L_1 \cup L_2;t) = \Phi(L_1;t) \cdot \Phi(L_2;t).$
\end{enumerate}

\medskip

This Jones polynomial satisfies a skein relation slightly different from
the usual one  and is normalized by dividing
by  $(t^{{1}/{2}}+ t^{-{1}/{2}})^{\mu}$  (see p.11, \cite{lin-wang}).
For example, if $L$ is Brunnian
(i.e. all proper sublinks are trivial)
and $J(L;t)$ is the usual Jones polynomial
(i.e. the one tabulated in C. Adams' book \cite{adams}),  then
$$\Phi(L;t) = (-1)^{\mu} \left( 1- \frac{J(L;1/t)}{(t^{1/2}+
t^{-1/2})^\mu}\right).$$
This and property (ii) is all we need to know to settle
the conjecture.

Define $$a_n(L)=\frac{1}{n!} \left. \frac{d^n \Phi(L;t)}{dt^n}\right|_{t=1}.$$
We claim that $a_n(L) \in \zz$ for $L$ a GSL.
For knots $K,$ it is evident that
$\Phi(K;t)$ is a Laurent polynomial, and property (ii) implies the same
for $L$ a GSL, i.e. $\Phi(L;t) \in \zz[1/t,t].$
Obviously the nth derivative of $t^m$ at $t=1$ is simply
$m(m-1)\cdots (m-n+1),$ which is divisible by $n!,$ since
any product of $n$ successive integers is divisible by $n!$
(p. 63 \cite{hardy-wright}).
It now follows by linearity that
$\left. \frac{d^n \Phi(L;t)}{dt^n}\right|_{t=1}$
is divisible by $n!,$ hence $a_n(L) \in \zz.$
This implies Conjecture 4.1 for GSLs,
but in fact, more is true.

\medskip

For any ASL $L$, Theorem 4.1 and Lemma 4.2 of \cite{lin-wang} show that
\begin{equation} \label{eq1}
\Phi(L;t)=\sum_{i=\mu+1}^\infty a_i(K) (t-1)^i,
\end{equation}
(i.e. $a_i(L)=0$ for $i \leq \mu$), and that
both $a_{\mu+1}(L)$ and $2 a_{\mu+2}(L) \in 3 \zz.$
But if $\mu=1$ and $L=K$ is a knot,
then $a_n(K) \in \zz,$ hence
$a_3(K) \in 3 \zz.$

\begin{prop} Suppose $L=K_1 \cup \cdots \cup K_\mu$ is a GSL. Then
$$\Phi(L;t)= \sum_{i=2 \mu}^\infty a_i(L) (t-1)^i,$$ where
$$\begin{array}{ll}
a_i(L) \in 3^\mu \zz& \hbox{ for } \; 2\mu \leq i \leq 3\mu, \hbox{ and} \\
a_i(L) \in 3^{4\mu-i} \zz& \hbox{ for } \;  3 \mu < i \leq 4 \mu,\\
a_i(L) \in  \zz& \hbox{ for } \;  4 \mu < i.
\end{array}$$
\end{prop}
\begin{pf} By property (ii), we have
$$\Phi(L;t) = \Phi(K_1;t) \cdots \Phi(K_\mu;t).$$
Multiplication of the series expansions of $\Phi(K_j;t) = \sum_{i=0}^\infty
a_i(K_j) (t-1)^i$ gives the formula
$$a_i(L)=\sum_{\sigma_1+\cdots+\sigma_\mu=i} a_{\sigma_1}(K_1) \cdots
a_{\sigma_\mu}(K_\mu).$$
The proposition now follows from the fact that
$a_0(K_j)=0=a_1(K_j),$ and that $a_2(K_j)$ and $a_3(K_j)$ are multiples of $3$
for $j=1,\ldots,\mu.$
\end{pf}

\medskip
\noindent
{\it Examples.} \quad (I) If $K$ is the left-hand trefoil, then
\begin{eqnarray*}
\Phi(K;t) &=& -t^4+t^3+t-1\\
&=&-3(t-1)^2-3(t-1)^3-1(t-1)^4.
\end{eqnarray*}
Taking $L$ to be the GSL $K \cup \cdots \cup K$ shows that
Proposition 1 is sharp.

\medskip
\noindent
(II) Let $L$ be the Whitehead link.
Consulting link tables\footnote{ $L=5^2_1$ in the standard notation
\cite{adams}. Note that $J(L;t)$
depends on a choice of orientation.},
we obtain
\begin{eqnarray*}
\Phi(L;t)&=&\frac{-t^{7/2}+2t^{5/2}-t^{3/2}+2t^{1/2}-t^{-1/2}+t^{-3/2}}{t^{1/2}+t^{-1/2}} -1\\
&=&-t^3+3t^2-4t+5+ t^{-1}- {8}(t+1)^{-1}. \\
&=&  \frac{-3(t-1)^3}{2}+\sum_{n=4}^\infty \frac{(-1)^n
(2^{n-2}-{1})(t-1)^n}{2^{n-2}}.
\end{eqnarray*}
In particular, $a_5(L)= -7/8$ and so $3! \, \phi_3(L)= 21,$ which
provides a counterexample to Conjecture 4.1.
Notice moreover that $\phi_7(L) = \frac{127}{32},$ thus
$n! \, \phi_n(L)$ need not be an integer.

\begin{prop}
If $L$ is an ASL with $\mu$ components, then
$$2^{n-2} a_n(L) \in \zz.$$
\end{prop}
\begin{pf}
For $L$ a GSL, this is a consequence of  Proposition 1, while for $L$ an ASL,
it follows
by induction on $n, \mu,$ and the double unlinking number, which is the number
of double crossings needed to change $L$
to a GSL, as we now explain.

Consider two crossings, one positive and
the other negative, between distinct components of $L.$
Write $L=L_{+-},$ and notice that  $L_{-+}$ is an ASL with $\mu$ components,
and that $L_{0+}$ and $L_{+0}$ are also ASL links, but with $\mu-1$ components.
Using the skein relation twice and subtracting,
we obtain the double crossing change formula
(cf. Lemma 4.1, \cite{lin-wang})
$$(t+1)\big(\Phi(L_{+-};t)-\Phi(L_{-+};t) \big) = (t^2-t)
\big(\Phi(L_{0+};t)-\Phi(L_{+0};t) \big).$$
Equating coefficients of the power series in equation (\ref{eq1})
gives
\begin{eqnarray*}
&  a_n(L_{+-}) & =  a_n(L_{-+}) + \frac{a_{n-1}(L_{-+}) - a_{n-1}(L_{+-})}{2}\\
& +& \frac{a_{n-1}(L_{0+}) - a_{n-1}(L_{+0})}{2} + \frac{a_{n-2}(L_{0+}) -
a_{n-2}(L_{+0})}{2}.
\end{eqnarray*}
The proposition now follows from this formula by induction, since
we can assume that it has already been established for
$L_{-+}, L_{0+},$ and $L_{+0},$ and that
$2^{n-3} a_{n-1}(L_{+-}) \in \zz.$
Note that the previous example indicates that
this result is sharp.
\end{pf}


\begin{thebibliography}{99}

\bibitem[A]{adams} C. Adams,
	{The Knot Book,}
	(1994) W. H. Freeman and Co., New York.
\bibitem[HW]{hardy-wright} G. H. Hardy and E. M. Wright,
	{The Theory of Numbers,}
	(1968) fourth edition, Oxford University Press, London.
\bibitem[LW]{lin-wang}  X.-S. Lin and Z. Wang,
       	{\em On Ohtsuki's Invariants of integral homology 3-spheres, I,}
        (1995) preprint.
\bibitem[O]{ohtsuki}  T. Ohtsuki,
        {\em A polynomial invariant of integral homology 3-spheres,}
        Math. Proc. Camb. Phil. Soc., (1995), {\bf 117}.

\end{thebibliography}
\end{document}